\documentclass{amsart}
\usepackage{amssymb}
\usepackage{amsmath}
  \usepackage{paralist}
  \usepackage{mathrsfs}
 \usepackage[colorlinks=true]{hyperref}
  \hypersetup{urlcolor=blue, citecolor=red}

\newtheorem{theorem}{Theorem}[section]

\newtheorem{proposition}{Proposition}[section]

\theoremstyle{definition}

\def\O{\Omega}

\def\F{{\bf{F}}}
\def\n{\nu}

\def\W0{W_0^{1,p}(\O)}

\def\va{\rightarrow}

\def\R{\mathbb R}

\def\rn{{{\R}^n}}
\def\rnpiu{{{\R}^n_+}}

\title[Sharp Hardy inequalities in the half space with trace remainder term]{Sharp Hardy inequalities in the half space\\
 with trace remainder term}

\author[A. Alvino - A. Ferone - R. Volpicelli]{}


 \email{alvinoan@unina.it}
 \email{adele.ferone@unina2.it}
 \email{roberta.volpicelli@unina.it}

\begin{document}
 \maketitle

 \centerline{\scshape Angelo Alvino and Roberta Volpicelli}
\medskip
{\footnotesize
 \centerline{Dipartimento di Matematica e Applicazioni \lq\lq R. Caccioppoli"}
   \centerline{Universit\`a di Napoli}
   \centerline{ Complesso Monte S. Angelo, Via Cintia Napoli - ITALY}
} 
\bigskip\bigskip
\centerline{\scshape Adele  Ferone }
\medskip
{\footnotesize
 \centerline{ Dipartimento di Matematica}
   \centerline{Seconda Universit\`a di Napoli}
   \centerline{Viale Lincoln 5,  81100 Caserta,  ITALY}
}

\date{}

\bigskip

\begin{abstract}
In this paper we deal with a class of inequalities which interpolate the Kato's inequality and the Hardy's inequality in the half space. Starting from the classical Hardy's inequality in the half space $\rnpiu =\R^{n-1}\times(0,\infty )$, we show that, if we replace the optimal constant $\frac{(n-2)^2}{4}$ with a smaller one $\frac{(\beta-2)^2}{4}$, $2\le \beta <n$, then we can add an extra trace-term equals to that one that appears in the Kato's inequality. The constant in the trace remainder term is optimal and it tends to zero when $\beta$ goes to $n$, while it is equal  to  the optimal constant in the Kato's inequality when $\beta=2$.
\end{abstract}
\section{Introduction}

Sobolev spaces play a fundamental rule in the study of differential and integral operators especially for their imbedding characteristics. Most of the embedding results  assert that, if $\Omega$ is an open set of $\rn$ with smooth boundary, then  $W^{1,p}(\Omega )$ is imbedded into some Lebesgue spaces $L^q(\Omega )$ or $L^q(\partial \Omega )$, with $q>p$. Here we want to consider the particular case $p=2$ and $\Omega=\rnpiu = \{(x,t)\in \rn : x\in \R^{n-1}, t>0\}$, $n\ge 3$, that is the upper half $n$-dimensional euclidean space.

The standard trace embedding theorem  asserts that if $u$ is any real valued function on $\rnpiu$, sufficiently smooth up to the boundary and decaying fast enough at infinity, then the trace of $u$ belongs to the Lebesgue space $L^{2^\star}(\partial \rnpiu)$, with $2^\star =\frac{2(n-1)}{n-2}$ (cf., e.g., \cite{A}). More precisely,  the following inequality holds:
\begin{equation}\label{esc}
\left ( \frac{n-2}{2}\right )^{1/2}\omega_{n-1}^{1/{2(n-1)}}\| u\|_{L^{2^\star}(\partial \rnpiu)}\le \|\nabla u\|_{L^2( \rnpiu)}\quad \forall u\in W^{1,2}(\rnpiu)
\end{equation}
where $\omega_{m}$ is, now and in the following,  the Lebesgue measure of the unit ball in $\R^{m}$. The constant that appears  in (\ref{esc}) together with the extremal functions
\begin{equation}\label{esc2}
u_a(x,t)=\left [ \left (a+t\right )^2+|x|^2\right ]^{-\frac{n}{2}+1}\quad a\in (0,\infty ).
\end{equation}
were found in 1988 by Escobar (\cite{E}, see \cite{Na} for the case $p\neq 2$).

Although optimal in  the framework of Lebesgue spaces the trace imbedding (\ref {esc}) admits an improvement in terms of Lorentz spaces. Indeed the following Kato's inequality holds
 \begin{equation}\label{kato}
2\frac{\Gamma^2 \left ( \frac{n}{4} \right )} {\Gamma^2 \left ( \frac{n-2}{4}\right )} \int_{\partial\rnpiu}\frac{u^2(x,0)}{|x|}dx\le \int_{\rnpiu} |\nabla u|^2(x,t)dxdt\qquad u\in W^{1,2}(\rnpiu )
\end{equation}
where, here and in the following, $\Gamma$ is the usual Gamma function defined as $\Gamma (s)=\int_0^\infty t^{s-1}\exp (-t)dt$.
The constant given in (\ref{kato}) is optimal but, unlike what happens for the standard trace inequality (\ref{esc}), it is never attained (\cite{He}, \cite{DDM}). Indeed the functions that are candidates to be extremals, are proportional to the solution of the problem
\begin{equation}\label{P0}
\begin{cases}
\Delta \varphi=0 & \textrm {in}\> \>\rnpiu,\\
\varphi  =|x|^{-\frac{n}{2}+1}& \textrm {on}\>\> \partial \rnpiu.\\
\end{cases}
\end{equation}
(then do not belong to $W^{1,2}(\rnpiu)$) and can be expressed   in terms of Legendre functions (see \cite{DDM}). Inequality (\ref{kato}) is an improvement  of (\ref{esc}), in the following sense: from inequality (\ref{esc}) we deduce that  if $u\in W^{1,2}(\rnpiu)$ then the trace of $u$
belongs to the Sobolev space $L^{2^\star}(\partial \rnpiu )$. In fact, inequality (\ref{kato}) tells us that the trace of $u$ has a higher degree of summability, since one can deduce that it  belongs to the Lorentz space $L^{2^\star, 2}(\partial \rnpiu)$ which is a proper subspace of $L^{2^\star}(\partial \rnpiu )$.

The relation between (\ref{esc}) and (\ref{kato}) is exactly the same that there exists between the classical Sobolev inequality
\begin{equation}\label{sob}
\sqrt{\pi n(n-2)}\left (\frac{\Gamma \left (\frac{n}{2}\right )}{\Gamma (n)}\right )^{1/n} \| u \|_{L^{\frac{2n}{n-2}}(\rn)} \le \| \nabla u \|_{L^2(\rn )}
\end{equation}
and the Hardy inequality
\begin{equation}\label{har}
\frac {\left (n-2\right )^2}{4}\int_\rn \frac{u^2(y)}{|y|^2}dy \le \int_\rn |\nabla u|^2(y)dy.
\end{equation}
Indeed, the constants that appear in (\ref{sob}) and (\ref{har}) are optimal (\cite{Au}, \cite{Da}, \cite{Ta}, see also  \cite{Al1}, \cite{CNV}, \cite{OK}), but only that one in (\ref{sob}) is attained.           Moreover, inequality (\ref{har}) strengthens the standard Sobolev embedding of $W^{1,2}(\rn )$ into $L^{2^*}(\rn)$, with $2^*=\frac{2n}{n-2}$ in as much as it states that $u\in L^{2^*,2}(\rn)$ and $ L^{2^*,2}(\rn)\subsetneq L^{2^*}(\rn)$.

Both the theories of boundary trace in Sobolev spaces and of Hardy inequalities have a large number of applications, especially to boundary value problems for partial differential equations and non linear analysis. They have been developed, via different methods and in different settings, by various authors ( see,  for example, \cite{Ad}, \cite{FGR}, \cite{FS},  \cite{GS} or the references on this topic in the monographs \cite{AH}, \cite{Ma}, \cite{OK}, \cite{Zi}). In particular, the lack of extremals in (\ref{kato}) and in (\ref{har}) has inspired  many mathematicians to consider possible extra terms on the left hand side.   As regards inequality (\ref{har}), it has been proven that no extra terms can be added on the left hand side (\cite{CF1}, \cite{CF2}), while, if $\rn$ is replaced by a bounded open subset containing the origin, then different type of remainder terms can be considered  (see, for example, \cite{ACR, AE, AVV, BFT, BWW, BM, BMS, BV, DHA, DDM, DELV, FMT1, FMT2, GGM, GGS, GM, HHL, Ti, VZ}). A similar result for inequality (\ref{kato}) has been proven in \cite{DDM} where a Kato's inequality with a remainder term is proven in the intersection of $\rnpiu$ with a ball centered at the origin.

 In this paper we deal with a class of inequalities which interpolate inequalities (\ref{kato}) and (\ref{har}). Starting from the classical Hardy's inequality in the half-space, we show that, if we replace the optimal constant with a smaller one, then we can add an extra term equal to that one that appears on the left hand side of (\ref{kato}).  Indeed, our aim is to prove   the following Theorem.
\begin{theorem}\label{main}
Let $n\ge 3$ and let $u$ be a real function on $\rnpiu$ vanishing at infinity, such that $|\nabla u|\in L^2 (\rnpiu)$. Then, for any $2\le \beta < n$, there exists a positive constant $H(n,\beta )$ such that
\begin{equation}\label{Hardy}
 H(n,\beta)
 \int _{\partial {\rnpiu}}\frac {u^2}{|x|} dx+\frac {(\beta-2)^2}{4}\int_ {{\rnpiu}}\frac {u^2}{|x|^2+t^2} dx dt\le  \int_ {{\rnpiu}}\vert \nabla u\vert^2 dx dt.
\end{equation}
The best value of the  constant $H(n,\beta)$ is given by
\begin {equation}\label{constant}
H(n,\beta)=2\frac{\Gamma \left ( \frac{n+\beta}{4}-\frac{1}{2}\right )\Gamma \left ( \frac{n-\beta}{4}+\frac{1}{2}\right )} {\Gamma \left ( \frac{n+\beta}{4}-1\right )\Gamma \left ( \frac{n-\beta}{4}\right )}. \end{equation}
\end{theorem}

\noindent It can be easily checked that when $\beta=2$, then inequality (\ref{Hardy}) reduces to (\ref{kato}), while when $\beta$ goes to $n$, then inequality (\ref{Hardy}) reduces to (\ref{har}).
The optimal constant $H(n,\beta )$, as  expected, is  never attained since the candidates to be extremal  functions are proportional to the solution of  the problems
\begin{equation}\label{P}
\begin{cases}
\Delta \varphi+\frac {(\beta-2)^2}{4}\frac  {\varphi}{|x|^2+t^2}=0 &\textrm {in} \>\>\rnpiu,\\
\varphi  =|x|^{-\frac{n}{2}+1}& \textrm {on}\>\> \partial \rnpiu.\\
\end{cases}
\end{equation}
These   solutions are explicitly given in Section 2 and they  are expressed in terms of the hyper-geometric series
\begin{equation} \label{series}
F(a,b,c;z)=1+\frac{\Gamma (c)}{\Gamma (a)\Gamma (b)}\sum_{k=1}^\infty \frac{ \Gamma (a+k)\Gamma (b+k)}{\Gamma (c+k)}\frac { z^k}{ k!}.
\end{equation}
The following Proposition holds.
\begin{proposition}\label{extremal}
Let $2\le \beta < n$ and let $H(n,\beta )$ be the constant defined in (\ref{constant}). Then the functions
\begin{equation}\label{extr}
\varphi (x,t) = \frac{F\left( \frac{n+\beta}{4}-1,\frac{n-\beta}{4},\frac{1}{2};\frac{t^2}{|x|^2+t^2} \right )}{\left [ |x|^2+t^2\right ]^{\frac{n-2}{4}}}  - t\frac{H(n,\beta )}{\left [ |x|^2+t^2\right ]^{\frac{n}{4}}} F\left(\frac{n+\beta}{4}-\frac {1}{2},\frac{n-\beta}{4}+\frac {1}{2},\frac{3}{2};\frac{t^2}{|x|^2+t^2}\right)
\end{equation}
are regular solutions of problems (\ref{P}).
 \end{proposition}

Incidentally, when  $\beta =2$ a much more handle expression of the solution of (\ref{P}) can be considered. Indeed     such  solution is proportional to the harmonic function
 \begin{equation}\label{phi}
 \phi(x,t)=\int_0^\infty \frac{a^{\frac{n}{2}-2}}{\left [ \left (a+t\right )^2+|x|^2\right ]^{\frac{n}{2}-1}}da.
 \end{equation}
 obtained by integrating the harmonic functions given in  (\ref{esc2}) with respect to a suitable weight of the parameter. The choice of the weight is influenced by the required summability on the boundary of $\rnpiu$.

 To conclude, let us describe our approach, based on a very classical method of Calculus of Variation (see \cite{Sa}, p. 167) that recently have been adopted in (\cite{Al2}) to find an improvement of classical Sobolev inequality. This method can be used  to prove many different integral inequalities and differs from that one used by Herbst,  based on dilation analytic techniques, and from that one used in \cite{DDM}, based on a suitable change of variables (see also \cite{BM,BV,GGM,Ma}). Our approach, instead, is based on the careful construction of a suitable divergence free vector field $\bf F$ related to the family $\mathcal G=\{k\varphi \}_{k\in(0,\infty )}$, where $\varphi$ is defined in (\ref{extr}). The functions of $\mathcal S$ are solutions of the Euler equation of the functional
\begin{equation}\label{J}
J(u)= \int_ {{\rnpiu}}\vert \nabla u\vert^2 dx dt-\frac {(\beta-2)^2}{4} \int_ {{\rnpiu}}\frac {u^2}{|x|^2+t^2} dx dt
\end{equation}
but, unfortunately,
they do not belong to $W^{1,2}(\rnpiu )$, since they
do not  have   the right summability neither at the origin nor at infinity. To overlap this inconvenient, we have to  evaluate $J(u)$  on an approximating sequence of bounded sets that do not contain the origin. The very definition of $\F$ allows us to   estimate $J(u)$ from below with the flow of $\F$ across the graph of $u$.
Since $\bf F$ is divergence free, when $u$ is sufficiently smooth, we can use divergence theorem to prove that the inequality (\ref{Hardy}) holds. The optimality of the constant is derived by repeating the previous arguments,  replacing  the  function $u$ by the functions  of $\mathcal G$ and observing that, in such case, all the inequalities hold as equalities.

\section{Proof of Proposition \ref{extremal}}

Writing
\begin{equation}\label{change}
\rho=\sqrt {|x|^2+t^2},\qquad \theta={\rm arctan}\frac  {t}{|x|},\qquad \varphi(x,t)=\rho^{-\frac {n}{2}+1}f(\theta)
\end{equation}
problem  (\ref {P})  is equivalent to  the  following limit problem
\begin{equation}\label{equ. f}
\begin{cases}
f''(\theta)-(n-2)\tan \theta f'(\theta)-\left(\frac {(n-2)^2}{4}-\frac {(\beta-2)^2}{4}\right)f(\theta)=0&\theta\in(0,\frac{\pi}{2})\\
f(0)=1 \quad \displaystyle \lim_{\theta\rightarrow \frac{\pi}{2}}f(\theta)\in\R.\\
\end{cases}
\end{equation}
 Equation (\ref {equ. f}) is explicitly  solved in \cite{PZ} (p. 271 eq. 131). Indeed  $f(\theta)=w(\sin^2\theta)$,  where $w$ is the solution of following limit problem for the the hypergeometric equation  (see, for instance, \cite{AS} ,\cite{PZ} for the general theory)
\begin{equation}\label{hyper}
\begin{cases}
z(z-1)w''(z)+\left(\frac{n}{2}z-\frac {1}{2}\right)w'(z)+\left[\frac {(n-2)^2}{16}-\frac {(\beta-2)^2}{16}\right]w(z)=0 & z\in(0,1)\\
w(0)=1 \quad \displaystyle  \lim_{z\rightarrow 1}w(z)\in\R&\\
\end{cases}
\end{equation}
The  general solution of (\ref {hyper})  satisfying $w(0)=1$  has the form
\begin{equation}\label {solution}w(z)=F\left(\frac{n+\beta}{4}-1,\frac{n-\beta}{4},\frac{1}{2};z\right)+\alpha\sqrt z F\left(\frac{n+\beta}{4}-\frac {1}{2},\frac{n-\beta}{4}+\frac {1}{2},\frac{3}{2};z\right)\qquad \end{equation}
where $F(a,b,c;z)$ is the hypergeometric series given in (\ref{series})
which, ${\it a\> fortiori}$,  is convergent for $0\le z<1$.
Since we are looking for bounded solution of $(\ref{equ. f})$, we have to analyze the behavior of a  hypergeometric function near the point  $z=1$. To this aim consider  that  (see \cite {AS} p. 559)
\begin {equation}\label {formula1}F(a,b,c;1)=\frac {\Gamma (c)\Gamma(c-a-b)}{\Gamma (c-a)\Gamma(c-b)}\qquad if\qquad c-a-b>0;\end {equation}
\begin {equation}\label {formula2}\lim_{z\rightarrow 1} \frac {F(a,b,c;z)}{\ln(1-z)}=-\frac {\Gamma(a+b)}{\Gamma (a)\Gamma(b)}\qquad if\qquad  c-a-b=0;
\end {equation}
\begin {equation}\label {formula3}\lim_{z\rightarrow 1} \frac {F(a,b,c;z)}{(1-z)^{c-a-b}}=\frac {\Gamma (c)\Gamma(a+b-c)}{\Gamma (a)\Gamma(b)}\qquad if\qquad  c-a-b<0.
\end {equation}
An easy calculation shows that for  both the hypergeometric functions appearing  in (\ref{solution})  $c-a-b=\frac{3-n}{2}\le 0$. Let us,  now,   first examine the case $n=3$. Since
\begin{equation}\label {limite w}
\lim_{z\rightarrow 1} w(z)=\lim_{z\rightarrow 1} \ln (1-z)\left[\frac {F\left(\frac{\beta-1}{4},\frac{3-\beta}{4},\frac{1}{2};z\right)}{\ln(1-z)}+\alpha \sqrt z \frac {F\left(\frac{1+\beta}{4},\frac{5-\beta}{4},\frac{3}{2};z\right)}{\ln(1-z)} \right],
\end{equation}
using (\ref{formula2}) and formula $5.3.10$ p. $559$ of \cite {AS}, we get that the limit is finite if and only if
$$\alpha=- \frac {\Gamma\left  (\frac{1}{2}\right )\Gamma \left ( \frac{\beta+1}{4}\right )\Gamma \left ( \frac{5-\beta}{4}\right )}{\Gamma \left ( \frac{3}{2}\right )\Gamma \left ( \frac{\beta-1}{4}\right )\Gamma \left ( \frac{3-\beta}{4}\right )}=-H(3,\beta)
$$
where in the last equality we use the fact that $\Gamma (\frac{1}{2})=2 \Gamma \left (\frac{3}{2}\right )$ and the very definition of $H(n,\beta)$ given in (\ref{constant}).

\noindent In a similar way, if $n>3$ consider
\begin{equation}\label {limite 2}
\lim_{z\rightarrow 1} w(z)=\lim_{z\rightarrow 1} (1-z)^{\frac {3-n}{2}}
\left [ \frac{F\left(\frac{n+\beta}{4}-1,\frac{n-\beta}{4},\frac{1}{2};z\right)}{(1-z)^{\frac {3-n}{2}}}+\alpha \sqrt z \frac{F\left(\frac{n+\beta}{4}-\frac {1}{2},\frac{n-\beta}{4}+\frac {1}{2},\frac{3}{2};z\right)}{(1-z)^{\frac {3-n}{2}}}\right ].
\end{equation}
Using  (\ref{formula3}) instead of (\ref{formula2}), de l'H{o}pital theorem and differentiation formulas for hypergeometric functions  we get that  the limit is  finite if and only if
$\alpha=-H(n,\beta)$.

Finally,  if we choose  $\alpha=-H(n,\beta)$ in (\ref {solution}) and we take  into account (\ref{change}), we deduce the thesis since
\begin{equation}\label{estremale}
\varphi (x,t)=\frac{1}{\left( {|x|^2+t^2}\right)^{\frac{n-2}{4}}}w\left(\frac {t^2}{|x|^2+t^2}\right ).
\end{equation}


%

\section{Proof of Theorem \ref{main}}\label{sec3}
Let  $\varphi$ be the function defined in (\ref{estremale}):
since $\varphi$ is a solution of  problem (\ref{P}), it is a solution of the Euler Lagrange equation associated to the functional (\ref{J}). Consider the one-parameter family of ${\mathcal H}^n$-surfaces $\mathcal{G}=\left \{ G_k \right \}_{k\ge 0}$ given by the graphs of the functions $\varphi_k=k\varphi$. For each point $(x,t,v)\in\rnpiu\times\R_+$ there exist a unique $G_{\bar k}$ containing it, that one corresponding to $\bar k=\frac{v}{\varphi (x,t)}$. Therefore, to each point $(x,t,v)\in\rnpiu\times\R_+$ we can associate the Mayer field
\begin{equation}\label{p}
\left ( 1, {\bf p}(x,t,v)\right )\equiv \left ( 1, \frac{v}{\varphi (x,t)}\nabla \varphi(x,t)\right ).
\end{equation}
Since $\varphi$ is a solution of problem (\ref{P}), then it is easy to check that the vector field
\begin{equation}\label{F}
{\bf F}(x,t,v)\equiv\left ( 2\frac{v}{\varphi(x,t)}\nabla \varphi (x,t), \frac{v^2}{\varphi^2(x,t)}|\nabla \varphi|^2(x,t)+\frac {(\beta-2)^2}{4} \frac {v^2}{|x|^2+t^2}\right )
\end{equation}
is divergence free.
Let us now take a nonnegative function $u\in C_0^\infty (\rn)$ vanishing outside the ball $B_R({ 0})$ centered at the origin and of radius $R$. If $0<r<R$, it is easy to check that
\begin{multline}\label{3}
\int_{\rnpiu\setminus B_r({ 0})} |{\bf\nabla} u|^2(x,t)dxdt-  \frac {(\beta-2)^2}{4} \int_{\rnpiu\setminus B_r({ 0})}  \frac {u^2(x,t)}{|x|^2+t^2}dxdt\\
=  \int_{\rnpiu\setminus B_r({\bf 0})}\left [2<{\bf p},{\bf\nabla} u>_n-|{\bf p}|^2 + |\bf{\nabla} u -{\bf p}|^2 \right]-\frac {(\beta-2)^2}{4}\int_{\rnpiu\setminus B_r({\bf 0})} \frac {u^2(x,t)}{|x|^2+t^2}dxdt \\
  \ge \int_{\rnpiu\setminus B_r({ 0})}\left [2<{\bf p},{\bf\nabla} u>_n-|{\bf p}|^2\right]-\frac {(\beta-2)^2}{4}\int_{\rnpiu\setminus B_r({ 0})} \frac {u^2}{|x|^2+t^2}dxdt\\
=\int_{\rnpiu\setminus B_r({ 0})} <{\bf F}(x,t,u(x,t)), ({\bf\nabla} u(x,t), -1)>_{n+1}dxdt
\end{multline}
where ${\bf p} \equiv {\bf p} (x,t,u(x,t))$ is defined in (\ref{p}) and $<\cdot , \cdot >_m $  stands for the standard inner products in $\R^m$. The last integral in (\ref{3}) is the inward flow of $\bf F$ across the graph of $u$ on ${\rnpiu\setminus B_r({ 0})}$.
 Since $\bf F$ is divergence free and $\F (x,t,{ 0})\equiv 0$, by divergence theorem it follows that  this flow equals the sum of the outward flows across the two manifolds
 \begin{equation}\label{s1}
\Sigma_1=\left \{(x,t,v)\in\rnpiu\times\R : |x|> r,  \quad  t=0,  \quad 0\le  v\le u(x,0)\right \}
\end{equation}
\begin{equation}\label{s2}
\Sigma_2=\left \{(x,t,v)\in\rnpiu\times\R : |x|^2+t^2 =r^2, \quad 0\le v\le u(x,t)\right \}.
\end{equation}
Let us begin evaluating  the flow of $\bf F$ across $\Sigma_1$: since the only non zero component of the unit outward normal ${\bf\nu}$ to $\Sigma_1$ is  the $n$-th $\n_n$  and $\n_n=-1$, then
\begin{multline}\label{flusso1}
\int_{\Sigma_1}<\F, \nu >_{n+1}\>d{\mathcal H}^n= -2\int_{|x|>r} dx\int_0^{u(x,0)} v\frac {\varphi_t(x,0)} {\varphi(x,0)}dv
\\=-\int_{|x|>r} {u^2(x,0) }\frac {\varphi_t(x,0)} {\varphi(x,0)}dx
=H(n,\beta) \int_{|x|>r}\frac  {u^2(x,0) }{|x|}dx.
\end{multline}
Indeed the first equality follows from the very definition of ${\bf F}$ given by (\ref{F}). The  last one relies on  the fact that, by (\ref{P}), $\varphi(x,0)=|x|^{-\frac{n}{2}+1}$ and on an easy  computation that shows
$\varphi_t(x,0)=-H(n,\beta ){|x|^{-\frac{n }{2}}}$.

\noindent As regards the flow of $\F$ across $\Sigma_2$, since its outer unit normal  is ${\bf \nu}\equiv \left (-\frac{x}{\sqrt{|x|^2+t^2}}, -\frac{t}{\sqrt{|x|^2+t^2}}, 0\right )$, then
\begin{equation}\label{flusso2}
\left |\int_{\Sigma_2}<\F , {\bf \nu}>_{n+1}d{\mathcal H}^n\right |\le \int_{\Sigma_2}v\left |\frac{\varphi_\rho (x,t)}{\varphi (x,t)}\right |d{\mathcal H}^n=\frac{n-2}{2r}\int_{\Sigma_2}v d{\mathcal H}^n \le \frac{n-2}{8r}n\omega_nr^{n-1}\displaystyle{\sup_{\partial B^+_r({0})}}u^2(x,t)
\end{equation}
where $\varphi_\rho$ denotes the derivative of $\varphi$ in the radial direction. The equality is a trivial consequence of the definitions of $\F$ and $\varphi$, given in (\ref{F}) and (\ref{change}) respectively, and the last inequality follows from the facts that  $v\le u(x,t)$ on $\Sigma_2$ and  $u$ is bounded on the boundary of $B^+_r({ 0})=\rnpiu\cap B_r(0)$.  Collecting (\ref{3}), (\ref{flusso1}) and (\ref{flusso2}) we deduce that for some positive constant $C$, depending only on $n$,  the following inequality holds
\begin{equation*}
\int_{\rnpiu\setminus B_r({\bf 0})} |\nabla u|^2(x,t)dxdt -\frac {(\beta-2)^2}{4}\int_{\rnpiu\setminus B_r({ 0})}  \frac {u^2(x,t)}{|x|^2+t^2}dxdt \ge H(n,\beta)\int _{|x|>r} \frac {u^2}{|x|} dx + Cr^{n-2}
\end{equation*}
from which (\ref{Hardy}) follows directly on letting $r$ go to zero.

It remains to show that the constant that appears in (\ref{constant}) is optimal. To do this, let  $0<r<R$ and apply the previous arguments replacing $u$ by $\varphi$. Starting from (\ref{3}) we get
\begin{multline}\label{opt}
\int_{B_R^+(0)\setminus B_r({\bf 0})} |\nabla \varphi|^2(x,t)dxdt-\frac {(\beta-2)^2}{4}\int_{B_R^+(0)\setminus  B_r({ 0})}  \frac {\varphi^2(x,t)}{|x|^2+t^2}dxdt \\=H(n,\beta)\int _{r<|x|<R}\frac {\varphi^2(x,0)}{|x|} dx + \frac{n-2}{2r}\int_{{\mathcal S}_1}vd{\mathcal H}^n - \frac{n-2}{2R}\int_{{\mathcal S}_2}vd{\mathcal H}^n
\end{multline}
where, as before,  $B^+_R(0)=\rnpiu\cap B_R(0)$ and
\begin{equation*}
{\mathcal S}_1=\left \{ (x, t,v)\in\rnpiu\times\R :\> |x|^2+t^2=r^2,\quad 0\le v\le \varphi (x,t)\right\}
\end{equation*}
 \begin{equation*}
{\mathcal S}_2=\left \{ (x, t,v)\in\rnpiu\times\R :\> |x|^2+t^2=R^2,\quad 0\le v\le \varphi (x,t)\right\}.
\end{equation*}
It is easy to check that the last two integrals in (\ref{opt}) are equal. Indeed, by spherical coordinates and (\ref{estremale})
\begin{multline}\label{uni}
\frac{n-2}{2r}\int_{{\mathcal S}_1}\!vd{\mathcal H}^n =\frac{n-2}{2r}\int_{\partial B_1^+(0)}\!\!\!r^{n-1}\left ( \int_0^{\varphi(rx', rt')}\!\!\!vdv\right)d{\mathcal H}^{n-1}\\
=\frac{n-2}{4}r^{n-2}\int_{\partial B_1^+(0)}\varphi^2(rx',rt')d{\mathcal H}^{n-1} = \frac{n-2}{4}\int_{\partial B_1^+(0)}w^2\left ({t'^2}\right )d{\mathcal H}^{n-1}=\frac{n-2}{2R}\int_{{\mathcal S}_2}vd{\mathcal H}^n.
\end{multline}
Collecting (\ref{opt}) and (\ref{uni})  we deduce
\begin{equation*}
\displaystyle{\lim_{R\va\infty}\lim_{r\va 0}}\frac{\int_{B_R^+(0)\setminus B_r({\bf 0})} |\nabla \varphi|^2(x,t)dxdt}{\int _{r<|x|<R}\frac {\varphi^2(x,0)}{|x|} dx }=H(n,\beta)
\end{equation*}
that shows the optimality of the constant.



\begin{thebibliography}{99}

\bibitem {AS} (MR0167642)
 \newblock {\sc M. Abramowitz and I.A. Stegun},
 \newblock ``Handbook of mathematical functions with formulas, graphs, and mathematical tables",
 \newblock   U.S. Government Printing Office, Washington, D.C., 1964.

\bibitem {A} (MR0450957)
 \newblock {\sc R.A. Adams},
 \newblock ``Sobolev spaces",
 \newblock   Academic Press, New York, 1975.

 \bibitem {AH} (MR1411441)
 \newblock {\sc R.A. Adams and L.I. Hedberg},
 \newblock ``Function spaces and potential theory",
 \newblock   Springer, Berlin, 1976.


\bibitem{Ad} (MR1918752)
 \newblock {\sc Adimurthi,}
 \newblock\emph{Hardy-Sobolev inequality in $H^1(\Omega)$ and its applications},
  \newblock Commun. Contemp. Math.,  \textbf{4} (2002), no. 3, 409-434.

\bibitem {ACR} (MR1862130)
 \newblock {\sc Adimurthi, N. Chaudhuri and M. Ramaswamy}
 \newblock\emph{An improved Hardy-Sobolev inequality and its application},
 \newblock Proc. Amer. Math. Soc. \textbf{130} (2002),  no. 2, 485-505.

 \bibitem {AE} (MR2184082)
 \newblock {\sc Adimurthi and M.J. Esteban,}
 \newblock\emph{An improved Hardy-Sobolev inequality in $W^{1,p}$W1,p and its application to Schršdinger operators},
 \newblock NoDEA Nonlinear Differential Equations Appl. {\bf 12} (2005), no. 2, 243-263.

 \bibitem {Al1} (MR0438106)
 \newblock {\sc A. Alvino,}
 \newblock\emph{Sulla diseguaglianza di Sobolev in spazi di Lorentz},
 \newblock Boll. Un. Mat. Ital. A (5) {\bf 14} (1977), no. 1, 148-156.

 \bibitem {Al2} (MR2550852)
 \newblock {\sc A. Alvino,}
 \newblock\emph{On a Sobolev-type inequality},
 \newblock Atti Accad. Naz. Lincei Cl. Sci. Fis. Mat. Natur. Rend. Lincei (9) Mat. Appl. {\bf 20} (2009), no. 4, 379-386.

 \bibitem {AVV} (MR2738657)
 \newblock {\sc A. Alvino, R. Volpicelli and B. Volzone}
 \newblock\emph{On Hardy inequalities with a remainder term},
 \newblock Ric. Mat. {\bf 59} (2010), no. 2, 265-280.

 \bibitem{Au} (MR0448404)
 \newblock {\sc T. Aubin,}
 \newblock\emph{Problmes isopŽrimŽtriques et espaces de Sobolev},
  \newblock J. Differential Geometry.,  \textbf{11} (1976), no. 4, 573-598.

\bibitem{BFT} (MR1970026)
\newblock {\sc G. Barbatis, S. Filippas and A. Tertikas,}
\newblock\emph{Series expansion for $L^p$ Hardy inequalities},
\newblock Indiana Univ. Math. J. {\bf 52} (2003), no. 1, 171-190.

 \bibitem{BWW} (MR2018667)
\newblock {\sc T. Bartsch, T. Weth, Tobias and M. Willem,}
\newblock\emph{A Sobolev inequality with remainder term and critical equations on domains with topology for the polyharmonic operator},
\newblock Calc. Var. Partial Differential Equations {\bf 18} (2003), no. 3, 253-268.

  \bibitem{BM} (MR1655516)
 \newblock {\sc H. Brezis and M. Marcus,}
 \newblock\emph{Hardy's inequalities revisited},
  \newblock Ann. Scuola  Norm. Sup. Pisa Cl. Sci.(4) \textbf{25}(1997), no. 1-2, 217-237 (1998).

   \bibitem{BMS} (MR1742864)
 \newblock {\sc H. Brezis, M. Marcus and I. Shafrir,}
 \newblock\emph{Extremal functions for Hardy's inequality with weight},
  \newblock J. Funct. Anal. {\bf 171} (2000), no. 1, 177-191.

   \bibitem{BV} (MR1605678)
 \newblock {\sc H. Brezis and J.L. Vazquez,}
 \newblock\emph{Blow-up solutions of some nonlinear elliptic problems},
  \newblock Rev. Mat. Univ. Complut. Madrid {\bf 10} (1997), no. 2, 443-469.

     \bibitem{CF1} (MR2541359)
 \newblock {\sc A. Cianchi and A. Ferone,}
 \newblock\emph{Best remainder norms in Sobolev-Hardy inequalities},
  \newblock Indiana Univ. Math. J. {\bf 58} (2009), no. 3, 1051-1096.

     \bibitem{CF2}
 \newblock {\sc A. Cianchi and A. Ferone,}
 \newblock\emph{Improving sharp Sobolev type inequalities by optimal remainder gradient  norms},
  \newblock to appear on Commun. Pure Appl. Anal..

 \bibitem{CNV} (MR2032031)
 \newblock {\sc D. Cordero-Erasquin, B. Nazaret and C. Villani,}
 \newblock\emph{A mass-transportation approach to sharp Sobolev and Gagliardo-Nirenberg inequalities},
  \newblock Adv. Math.  \textbf{182}(2004), 307--332.

\bibitem{Da} (MR1747888 )
 \newblock {\sc E.B. Davies,}
 \newblock\emph{A review of Hardy inequalities},
  \newblock The Maz'ya anniversary collection, Vol. 2 (Rostock, 1998),
Oper. Theory Adv. Appl.,   Vol. 110, Birkhauser Basel, 1999, 55-67.

\bibitem{DDM} (MR2393398)
 \newblock {\sc J. Davila, L. Dupaigne and M. Montenegro,}
 \newblock\emph{The extremal solution of a boundary reaction problem},
  \newblock Commun.  Pure  Appl. Anal.,  \textbf{7} (2008), no. 4, 795-817.


 \bibitem{DHA} (MR2108051)
 \newblock {\sc A. Detalla, T. Horiuchi and H. Ando,}
 \newblock\emph{Missing terms in Hardy-Sobolev inequalities and its application},
  \newblock Far East J. Math. Sci. (FJMS) {\bf 14} (2004), no. 3, 333-359.

    \bibitem{DELV} (MR2091354)
 \newblock {\sc  J. Dalbeault, M.J. Esteban, M. Loss, and L. Vega}
 \newblock\emph{An analytical proof of Hardy-like inequalities related to the Dirac operator},
  \newblock J. Funct. Anal. {\bf 216} (2004), no. 1, 1-21.



\bibitem {E} (MR0962929)
 \newblock {\sc J.F. Escobar},
 \newblock\emph{Sharp constant in a Sobolev trace inequality,}
 \newblock  Indiana Univ. Math. J. \textbf{37} (1988), 687-698.

 \bibitem {FGR}  (MR2295124)
 \newblock {\sc J. Fernandez Bonder, P. Groisman and J.D. Rossi}
 \newblock\emph{Optimization of the first Steklov eigenvalue in domains with holes: a shape derivative approach,}
 \newblock  Ann. Mat. Pura Appl. \textbf{186} (2007), no.2,  341-358.

 \bibitem {FS} (MR2413374)
 \newblock {\sc J. Fernandez Bonder and N. Saintier,}
 \newblock\emph{Estimates for the Sobolev trace constant with critical exponent and applications,}
 \newblock  Ann. Mat. Pura Appl. \textbf{187} (2008), no.4,  683-704.

 \bibitem {FMT1} (MR2297247)
 \newblock {\sc S. Filippas, V.G. Maz'ya and A. Tertikas,}
 \newblock\emph{Critical Hardy-Sobolev inequalities,}
 \newblock  J. Math. Pures Appl. (9) {\bf 87} (2007), no. 1, 37-56.

  \bibitem {FMT2} (MR2214621)
 \newblock {\sc S. Filippas, V.G. Maz'ya and A. Tertikas,}
 \newblock\emph{On a question of Brezis and Marcus,}
 \newblock  Calc. Var. Partial Differential Equations {\bf 25} (2006), no. 4, 491-501.


\bibitem {GS} (MR2396523)
 \newblock {\sc M. Gazzini and E. Serra,}
 \newblock\emph{The Neumann problem for the HŽnon equation, trace inequalities and Steklov eigenvalues,}
 \newblock  Ann. Inst. H. PoincarŽ Anal. Non LinŽaire \textbf{25} (2008), no.2,  281-302.

\bibitem {GGM} (MR2048513)
 \newblock {\sc F. Gazzola, H.C. Grunau and E. Mitidieri,}
 \newblock\emph{Hardy inequalities with optimal constants and remainder terms,}
 \newblock  Trans. Amer. Math. Soc. {\bf 356} (2004), no. 6, 2149-2168.

 \bibitem {GGS} (MR2010961)
 \newblock {\sc F. Gazzola, H.C. Grunau and M. Squassina,}
 \newblock\emph{Existence and nonexistence results for critical growth biharmonic elliptic equations,}
 \newblock  Calc. Var. Partial Differential Equations {\bf 18} (2003), no. 2, 117-143.


 \bibitem {GM} (MR2443723)
 \newblock {\sc N. Ghoussoub and A.  Moradifam,}
 \newblock\emph{On the best possible remaining term in the Hardy inequality,}
 \newblock  Proc. Natl. Acad. Sci. USA {\bf 105} (2008), no. 37, 13746-13751.






\bibitem{He} (MR0436854)
\newblock{\sc I.W. Herbst,}
\newblock\emph{Spectral theory of the operator $(p^{2}+m^{2})^{1/2}-Ze^{2}/r$,}
 \newblock  Comm. Math. Phys. \textbf{53}  (1977), no. 3, 285-2940.

 \bibitem{HHL} (MR1892180)
\newblock{\sc M. Hoffmann-Ostenhof, T. Hoffmann-Ostenhof and A.  Laptev,}
\newblock\emph{A geometrical version of Hardy's inequality,}
 \newblock  J. Funct. Anal. {\bf 189} (2002), no. 2, 539-548.


\bibitem {Ma} (MR0817985)
\newblock {\sc V.J. Maz'ja,}
\newblock\emph{``Sobolev spaces"},
\newblock Transl. from the Russian by T. O. Shaposhnikova, Springer-Verlag, 1985.

\bibitem {Na} (MR2258478)
\newblock {\sc B. Nazaret}
\newblock\emph{Best constant in Sobolev trace inequalities on the half-space},
\newblock Nonlinear Anal. \textbf{65} (2006), 1977-1985.



\bibitem {OK} (MR1069756)
\newblock {\sc B. Opic and A. Kufner}
\newblock\emph{``Hardy-type inequalities"},
\newblock Longman Scientific \& Technical, Harlow, 1990.

\bibitem {PZ} (MR2001201)
\newblock {\sc A.D. Polyanin and V.F. Zaitsev}
\newblock\emph{``Handbook of exact solutions for ordinary differential equations"},
\newblock Second edition, Chapman and Hall/CRC, Boca Raton, FL, 2003.

\bibitem {Sa} (MR1210325)
\newblock {\sc H. Sagan}
\newblock\emph{``Introduction to the calculus of variations"},
\newblock Dover Publication, Inc., New York, (1992).

\bibitem {Ta} (MR0463908)
\newblock {\sc G. Talenti,}
\newblock\emph{Best Constant in Sobolev Inequality},
\newblock Ann. Mat. Pura Appl. \textbf{110} (1976),  353-372.

\bibitem {Ti} (MR2124873)
\newblock {\sc J. Tidblom,}
\newblock\emph{A Hardy inequality in the half-space},
\newblock J. Funct. Anal. {\bf 221} (2005), no. 2, 482-495.

\bibitem {VZ} (MR1760280)
\newblock {\sc J.L. Vazquez and E. Zuazua,}
\newblock\emph{The Hardy inequality and the asymptotic behaviour of the heat equation with an inverse-square potential},
\newblock J. Funct. Anal. {\bf 173} (2000), no. 1, 103-153.

\bibitem {Zi} (MR1014685)
\newblock {\sc W.P. Ziemer,}
\newblock\emph{``Weakly differentiable functions. Sobolev spaces and functions of bounded variation"},
\newblock Springer-Verlag, New York, 1989.

\end{thebibliography}
\end{document}